\input amstex.tex 
\documentstyle{amsppt} 
\magnification=1200 
\baselineskip 16pt plus 2pt 
\advance\vsize -2\baselineskip
\parskip 2pt 
\NoBlackBoxes

\topmatter 
 
\title Lectures on Witten Helffer Sj\"ostrand  theory \endtitle 
 
\thanks   Supported in part by NSF
\endthanks

\author D. Burghelea (Ohio State University)
\endauthor 
\leftheadtext {Witten Helffer Sj\"ostrand  theory} 
\rightheadtext {D. Burghelea} 
  
\abstract 

Witten- Helffer-Sj\"ostrand theory is a considerable addition to 
the De Rham- Hodge theory for Riemannian
manifolds and can serve as a general tool to prove results about comparison of
numerical invariants associated to compact manifolds analytically, 
i.e. by using a Riemannian metric, 
or combinatorially, i.e by using a triangulation. In this presentation  
a triangulation, or a partition of a smooth manifold in cells,
will be viewed in a more analytic spirit, being provided by the 
stable manifolds of the gradient of a nice Morse function.
WHS theory was recently used both for providing  new proofs for known but 
difficult results in topology, as well as 
new results and a positive solution for an important conjecture
about $L_2-$torsion, cf [BFKM]. This  presentation 
is a short version of a one quarter course I have given
during the spring of 1997 at OSU.

\endabstract 
\toc 
\widestnumber\subhead{3.1} 
\head 0. Introduction. \endhead 
\head 1. Triangulations from an analytic point of view. \endhead 
\head 2. De Rham theory and integration. \endhead 
\head 3. Witten deformation and the main results of WHS-theory.\endhead
\head 4. The ideas of the proof.\endhead
\head 5. Extensions and a survey of applications. \endhead 
\head 6. References. \endhead 
\endtoc

\endtopmatter 
 
\document  
\vfill \eject 
 
\proclaim{0. Introduction}\endproclaim 
 
Witten Helffer Sj\"ostrand theory, or abbreviated, WHS -theory,  is a substantial addition to 
the De Rham -Hodge theory cf~\cite{DR} and a powerful tool for comparing numerical invariants 
associated to compact manifolds analytically (i.e by using a Riemannian 
metric,) and combinatorially 
(i.e by using a triangulation), cf[BZ1], [BZ2], [BFKM], [BFK1],
[BFK2]. It states in a precise way the relationship 
between the De Rham complex of a
manifold and the 
cochain complex provided by a smooth triangulation when described with the 
help of a Riemannian metric and of a Morse function. While there are other results which relate these two complexes, WHS-theory provides a connection 
between these two complexes with scalar producs and permits to relate some of the spectral properities of the Riemannian Laplacians given by the metric
and the combinatorial 
Laplacians given by the triangulation. 

The intuition behind  the WHS -theory is provided  by
physics and consists in regarding a compact smooth  manifold equipped with a Riemannian metric 
and a Morse function as an interacting system of harmonic oscillators. This 
intuition was first noticed 
and exploited by E. Witten, cf[Wi], in order to provide 
a short "physicist's proof " of Morse inequalities, a rather simple but very 
useful result in topology.  

Helffer and 
Sj\"ostrand have completed Witten's 
 picture with their results on Schr\"odinger operators and 
have considerably strengthened
Witten's mathematical statements, cf~\cite {HS2}. The work of Helffer and 
Sj\"ostrand on the Witten theory 
can be substantially simplified  by using simple observations familiar 
to topologists,
cf~\cite {BZ2} and ~\cite {BFKM}. 
As presented in ~\cite{HS}, their work, although very appealing,
is not very accessible to topologists  because of a large amount of 
estimates and preliminary results about 
Schr\"odinger operators. It turns out that not all of them are necessary and 
the Witten Helffer Sj\"ostrand work, at least as needed 
by topologists, can be presented 
and explained in a selfcontained manner and on a reasonable number of pages,
cf. section 5, [BFKM]. This survey is a presentation of the WHS-theory with these 
simplifications, hopefully accessible to a 
graduate student in geometry 
and topology and in a way appropriate to topological applications.

The mathematics behind the WHS-theory is almost entirely based on the following two facts: the 
existence of a gap in the spectrum of the Witten Laplacians detected 
by elementary mini-max characterization of the spectrum of selfadjoint positive operators and simple estimates involving the equations 
of the harmonic oscillator.  The 
Witten Laplacians in the neighborhood of critical points in ``admissible coordinates'' are given by such equations.

Witten's ideas  presented below were used by Witten to provide a new proof of Morse type inequalities 
and holomorphic Morse inequalities. WHS-theory was also used to provide a new proof of the equality of 
analytic and Reidemeister torsion and the $L_2$-version of WHS-theory to  
provide a new proof of the equality of 
Novikov- Shubin invariants defined analytically and combinatorially,cf [BFKM].

The $L_2-$version of WHS theory turned out to be not only important but so 
far an unavoidable ingredient in 
the proof of the equality of $L_2$ -analytic and Reidemeister torsion and of generalizations of this result, cf [BFKM].

\vfill \eject 
 
\subhead {1.Triangulations from an analytic point of view}
\endsubhead 
 
Let $M^n$ be a compact closed smooth manifold of dimension $n$.
A generalized triangulation is provided by a pair $(h,g),$ $h:M\to \Bbb R$ 
a smooth function, $g$ a Riemannian metric so that :

C1.  For any critical point $x$ of $h$ there exists a coordinate chart in the 
neighborhood of $x$ so that
in these  coordinates $h$ is quadratic  and $g$ is Euclidean.

Precisely, for any $x$ critical point of $h,$ ($x \in Cr(h)$), 
there exists a coordinate chart 
$\varphi: (U,x)\to (D_{\epsilon},0),\  U$ an open neighborhood of $x$ in $M,\   
D_{\epsilon}$ an open disc of 
radius $\epsilon$ in $\Bbb R^n,\  \varphi$ a diffeomorphism with 
$\varphi(x)=0,$
so that :
$$(i): h\cdot \varphi^{-1}(x_1, x_2, \cdots, x_n)=c - 1/2( x_1^2+\cdots x_k^2) + 
1/2 ( x_{k+1}^2+\cdots x_n^2)$$
$$(ii): (\varphi^{-1})^*(g)\  \text{is given by}\  g_{ij}(x_1,x_2,\cdots,x_n)= \delta_{ij}$$
Coordinates so that (i) and (ii) hold are 
called admissible.

It follows that any critical points has a well defined index, the number k
of the negative squares in the expression (i), which is independent of the 
choice of  a coordinate system with respect to which $h$ has the form (i).

C2. h is self indexing, i.e. for any critical point $x\in Cr(h)$
$h(x)= \text{index}\ x.$ 

Consider the vector field $-grad_g(h)$ and for any
$y\in M,$
denote  by $\gamma_y(t), -\infty < t < \infty,$ the unique trajectory of 
$-grad_g(h)$ which satisfies
the condition $\gamma_y(0)=y.$
\newline For $x\in Cr(h)$ denote by $W_x^-$ resp. $W_x^+$ the sets 
$$W_x^{\pm}=\{y\in M| \lim_{t\to \pm\infty}\gamma_y(t)=x\}.$$
In view of (i), (ii) and of the theorem of existence, 
unicity and smooth 
dependence on the initial condition 
for the solutions of ordinary differential equations, $W_x^-$
resp. $W_x^+$ is a smooth submanifold diffeomorphic to $\Bbb R^k$
resp. to $\Bbb R^{n-k},$ with $k=$index$x.$ This can be verified easily based on the fact that:
$$\varphi(W_x^-\cap U_x)=\{(x_1,x_2,\cdots,x_n)\in D(\epsilon) |x_{k+1}=x_{k+2}= \cdots= x_n=0\},$$
and 
$$\varphi(W_x^+\cap U_x)=\{(x_1,x_2,\cdots,x_n)\in D(\epsilon) |x_{1}=x_{2}= \cdots= x_k=0\}.$$ 

Since $M$ is compact and C1 holds, the set $Cr(h)$ is finite and since $M$ 
is closed (i.e. compact and without boundary), $M= \bigcup_{x\in Cr(h)}W_x^-.$ As already observed each $W_x^-$ is a smooth submanifold diffeomorphic to 
$\Bbb R^k,$ $k=$index $x,$ i.e. an open cell. 

C3. The vector field $-grad_g h$ satisfies the Morse-Smale condition if for any $x,y \in Cr(h)$, $W_x^- $ and $ W_y^+$
 are transversal. 

C3 implies that $\Cal M(x,y):= W_x^- \cap W_y^+$ is a smooth manifold of dimension equal to 
 $\text{Index}\ x- \text{Index}\ y.$ $\Cal M(x,y)$
is equipped with the action $\mu: \Bbb R\times \Cal M(x,y) \to
 \Cal M(x,y),$
defined by $\mu(t, z) = \gamma _z(t).$  

If $\text{Index}\ x \leq \text{Index}\ y,$
 and $x\neq y,$ in view of the transversality requested by the Morse Smale
condition, 
$\Cal M(x,y) = \emptyset.$

If $x\neq y$ and $\Cal M(x,y) \neq \emptyset,$
the action $\mu$ is free and we denote the quotient $\Cal M(x,y)/\Bbb R $
by $\tilde{\Cal M}(x,y);$
$\tilde \Cal M(x,y)$ is a smooth manifold of dimension $\text{Index}\ x-
\text{Index}\ y -1,$
diffeomorphic to the submanifold  
$\Cal M(x,y) \cap h^{-1}(c),$ for any real number $c$ in the open
interval $(\text{Index}\ x,  \text{Index}\ y).$ The elements of $\tilde \Cal M(x,y)$ are the trajectories from ``$x$ to $y$'' and such an element will be  
denoted by $\gamma.$ 

If $x=y,$ then $W_x^- \cap W_x^+ = {x}.$

The condition C3 implies that the 
partition of $M$ into open cells is actually a smooth cell complex.  
To formulate this fact precisely we 
recall that an 
\newline $n-$dimensional  manifold $X$ with corners is a paracompact  Hausdorff space  equipped with a 
maximal smooth atlas with charts  $\varphi : U\to \varphi(U)\subseteq \Bbb R^n_+$
with $\Bbb R^n_+= \{(x_1, x_2,\cdots x_n) | x_i\geq 0\}.$  The collection of points 
of $X$ 
which correspond (by some and then by any chart) to points in $\Bbb R^n$ 
with exactly $k$ coordinates equal to zero is
a well defined subset of $X$ and it will be denoted by $X_k.$ It 
has a structure of a smooth 
$(n-k)-$dimensional manifold. $\partial X = X_1 \cup X_2 \cup \cdots X_n$
is a closed subset which is a topological manifold and $(X,\partial X)$
is a topological manifold with boundary $\partial X.$ A compact 
smooth manifold
with corners, $X,$ with interior diffeomorphic to the Euclidean space,
will be called a compact smooth cell.

For any string of critical points $x=y_0, y_1,\cdots ,y_k$ with 
$$\text{index}\ y_0 > \text{index}\ y_1>,\cdots,> \text{index}\ y_k,$$ 
consider the smooth manifold of dimension  $\text{index}\ y_0 -k,$
$$\tilde{\Cal M}(y_0,y_1)\times\cdots \tilde{\Cal M}(y_{k-1},y_k)\times W_{y_k}^-,$$ and the smooth 
map $$i_{y_0, y_1,\cdots ,y_k}: \tilde{\Cal M}(y_0,y_1)\times\cdots 
\times \tilde{\Cal M}(y_{k-1},y_k)\times W_{y_k}^-
\to M,$$
defined by $i_{y_0, y_1,\cdots ,y_k}(\gamma_1,\cdots, \gamma_k, y):= i_{y_k}(y)$, for 
$\gamma_i \in \tilde{\Cal M}(y_{i-1}, y_i)$ and $y\in W_{y_k}^-,$ with $i_x: W_x^- \to M$ the inclusion of $W_x^-$ in $M.$

\proclaim {Theorem1.1}
Let $\tau =(h,g)$  be a generalized triangulation. For any critical point $x\in Cr(h)$
the smooth manifold $W_x^-$ has a canonical compactification $\hat W_x^-$ to a compact manifold 
with corners and the inclusion $i_x$ has a smooth extension 
\newline $\hat i_x: \hat W_x^- \to M$ so that :

\noindent (a): $(\hat W_x^-)_k = \bigcup_{(x,y_1,\cdots, y_k)} 
\tilde{\Cal M}(x,y_1)\times \cdots\times
\tilde{\Cal M}(y_{k-1}, y_k)\times W_{y_k}^-,$

\noindent (b): the restriction of $\hat i_x$ to $\tilde{\Cal M}(x,y_1)\times \cdots \times
\tilde{\Cal M}(y_{k-1}, y_k)\times W_{y_k}^-$ is given by 
\newline $i_{x=y_0, y_1\cdots, y_k}.$

\endproclaim

This theorem was probably well known to experts before it was formulated by Floer in the framework
of $\infty-$dimensional Morse theory cf.~\cite {F}. In fact, a  weaker version 
of this theorem, e.g. Proposition 2
in ~\cite {L}, suffices to conclude that the linear maps 
${\text{Int}^q}$'s 
defined in section 2 provide a morphism of cochain complexes. This is the  
only fact one needs in order to formulate the WHS-theory. However, Theorem 1.1 is a 
statement worth to be known. As formulated Theorem 1.1 is 
proven in [AB].

The name of generalized triangulation for $\tau= (h,g)$ is justified 
by the fact that any simplicial 
smooth triangulation  can be obtained as a generalized triangulation, 
cf~\cite{Po}.
We also point out that given a selfindexing Morse function $h$ 
and a Riemannian
metric $g,$ one can perform arbitrary small $C^0-$ perturbations to $g,$ 
so that the pair consisting of $h$ and the perturbed metric is a 
generalized triangulation, cf [Sm].

Given a generalized triangulation $\tau= (h,g),$ and for any critical
point $x\in Cr(h)$ an orientation 
$\Cal O_x$ of $W_x^-,$ one can associate a cochain complex of vector spaces 
over the field $\Bbb K$ or real or complex numbers, $(C^*(M,\tau), 
\partial^*).$ The diferential $\partial^*$ depends on the choosen 
orientations $\Cal O_x.$
To describe this complex we introduce the incidence numbers
$$ I_q: Cr(h)_q\times Cr(h)_{q-1} \to \Bbb Z$$ defined as follows:
\newline If $\tilde{\Cal M} (x,y)=\emptyset,$ we put $I_q(x,y)= 0.$
\newline If $\tilde{\Cal M} (x,y) \ne \emptyset,$  for any 
$\gamma\in \tilde{\Cal M} (x,y),$ the set  
$\gamma\times W_y^-$ appears as an open set of the boundary 
$\partial\hat{W}_x^-$ and the orientation 
$\Cal O_x$ induces an orientation on it. 
If this is the same as the orientation $\Cal O_y,$ 
we set  $\epsilon(\gamma)= +1,$ otherwise we set 
$\epsilon(\gamma)= -1.$ Define $I_q(x,y)$ by
$$I_q(x,y)= \sum_{\gamma
\in \tilde{\Cal M} (x,y)} \epsilon(\gamma).$$
In the case $M$ is an oriented manifold, the orientation of 
$M$ and the orientation $\Cal O_x$ on $W_x^-$ induce
an orientation $\Cal O_x^+$ on the stable 
manifold $W_x^+.$ For any $c\in (\text{index}\ y, \text{index}\ x),$
$h^{-1}(c)$ carries a canonical orientation induced from the orientation of  
$M.$ 
One can check  that $I_q(x,y)$ is the intersection number of 
$W_x^-\cap h^{-1}(c)$ with $W_y^+\cap h^{-1}(c)$ inside $h^{-1}(c)$
and is also the incidence number of the open cells 
$W_x^-$ and $W_y^-$
in the $CW-$ complex structure provided by $\tau.$

Denote by $(C^*(M,\tau),\partial^*)$ the cochain 
complex of $\Bbb K-$vector spaces defined by 

\noindent (1)  $C^q(M,\tau):= Maps (Cr_q(h), \Bbb K)$

\noindent (2)  $\partial^{q-1}: C^{q-1}(M,\tau) \to C^q(M,\tau), 
\  \partial ^{q-1}f (x)=
\sum _{y\in Cr_{q-1}(h)} I_q(x,y)f(y),$ 
where $x\in Cr_q(h).$  

Since $C^q(M,\tau)$ is equipped with a canonical base provided by the 
maps  $E_x$ defined by $E_x(y)= \delta_{x,y},$
$x,y\in Cr_q(h),$
it carries a natural scalar product which makes  
$E_x,$ $x\in Cr_q(h)$ orthonormal.

\proclaim {Proposition 1.2}
For any $q, \  \partial^q\cdot \partial ^{q-1}=0.$ 
\endproclaim
A geometric proof of this Proposition follows from Theorem 1.1 (cf [F] or
[AB]), 
The reader can also derive it by noticing that 
$(C^*(M,\tau),\partial^*)$ as defined is nothing but 
the cochain complex associated to the $CW-$complex structure 
provided by $\tau.$

\vfill \eject 
 
\proclaim{2.  De Rham theory and Integration theory}\endproclaim 
 
 Let $M$ be a closed smooth  manifold and $\tau= (h,g) $ be a generalized triangulation. 
 Denote by $(\Omega^*(M), d^*)$ the De Rham complex of $M.$ This is a 
cochain complex whose component $\Omega^r(M)$ is the (Frechet) space of 
 smooth differential forms of degree $r$ and whose differential 
$d^r: \Omega^r(M) \to \Omega^{r+1}(M)$
is given by the exterior differential. Recall that Stokes theorem can be formulated as follows:

 \proclaim{Theorem 2.1} Let $P$ be a compact $r-$dimensional oriented 
smooth manifold with corners 
 and $f:P\to M$ be a smooth map. Denote by $\partial f: 
P_1\to M $ the restriction of $f$ to 
 the smooth oriented manifold $P_1$ ($P_1$ defined as above).
If $\omega \in \Omega^{r-1}(M)$
 is a smooth form then $\int_{P_1} \partial f^*(\omega) $ 
is convergent and
 $$\int_P f^*(d\omega)= \int_{P_1} \partial f^*(\omega).$$
 \endproclaim

Consider the linear map 
$Int^q: \Omega^q(M) \to C^q(M,\tau),$ with
\newline $C^q(M,\tau)= 
\text{Maps}(Cr(h)_q,\Cal K),$ $Cr_q(h)= \{x\in Cr(h)| \text{index}\ x=q\}$ 
defined by
 $$Int^q(\omega)(x)= \int _{\hat W_x^-} \omega,$$

The collection of the linear maps ${Int^q}$' s defines a morphism 
$$Int^*: (\Omega^*(M), d^*) \to (C^*(M,\tau),\partial^*)$$
of cochain complexes.

\proclaim {Theorem 2.2}(De Rham) $Int^*$  induces an 
isomorphism in cohomology.

\endproclaim

Theorem 3.2, one of the two main results of the WHS theory,whose proof will be sketched below, 
is a considerable strengthening of this theorem. 

\vfill \eject 

\proclaim{3.  Witten deformation and the main results of WHS-theory}  \endproclaim 

Let $M$ be a closed manifold and $h:M\to \Bbb R$ a smooth function. For $t>0$ we consider
the complex $(\Omega^*(M), d^*(t))$
with differential $d^q(t):\Omega^q(M)\to \Omega^{q+1}(M)$ 
given by  $d^q(t)= e^{-th} d e^{th}$ or equivalently

$$
\leqalignno{ d^q(t)(\omega)= d\omega +t dh\wedge \omega.   &&(3.1)\cr}
$$  

$d^*(t)$ is the unique differential in $\Omega^*(M)$ which makes the multiplication 
by the smooth function $e^{th}$
an isomorphism of cochain complexes
$$e^{th}: (\Omega^*(M), d^*(t)) \to (\Omega^*(M), d^*).$$

Recall that for any vector field $X$ on $M$ one defines the zero 
order differential operator, 
$\iota_X= \iota_X^*: \Omega^*(M) \to \Omega^{*-1}(M),$ by 
$$\leqalignno{\iota_X^q\omega (X_1, X_2,\cdots, X_{q-1}) := \omega (X, X_1,\cdots, X_{q-1}) 
&&(3.2)\cr}$$
and the first order differential operator
$L_X= L_X^*:\Omega^*(M) \to \Omega^*(M),$ the Lie derivative in the 
direction $X$,

by 
$$\leqalignno{ L_X^q:= d^{q-1}\cdot \iota_X^q + \iota_X^{q+1}\cdot d^q.
&&(3.3)\cr}$$
They satisfy the following identities:
$$\leqalignno{\iota_X(\omega_1\wedge \omega_2)= \iota_X(\omega_1)\wedge \omega_2 +
(-1)^{|\omega_1|}\omega_1 \wedge \iota_X(\omega_2) 
.&&(3.4)\cr}$$
for $\omega_1\in \Omega^{|\omega_1|}(M),$ and
$$
\leqalignno{ L_X(\omega_1\wedge \omega_2)= L_X(\omega_1)\wedge \omega_2 +
\omega_1 \wedge L_X(\omega_2). 
&&(3.5)\cr}$$

Given a Riemannian metric $g$ on the oriented manifold $M$ we have the zeroth 
order operator 
$R^q:\Omega^q(M) \to \Omega^{n-q}(M)$, known as the star-Hodge operator
which, with respect to an oriented orthonormal frame 
$e_1, e_2,\cdots, e_n$ in the cotangent space at $x,$ is given by 
$$\leqalignno{ R_x^q (e_{i_1}\wedge\cdots \wedge e_{i_q}) = \epsilon (i_1,\cdots,i_q) e_1\wedge \cdots\wedge
\hat e_{i_1}\wedge\cdots \wedge\hat e_{i_q}\wedge\cdots\wedge e_n 
,&&(3.6)\cr}$$
$ 1\leq i_1<i_2,\cdots,i_q\leq n ,$ with $\epsilon (i_1,i_2,\cdots,i_q)$ 
denoting the sign of the
permutation of $(1,2\cdots n)$ given by 
$$(i_1,i_2,\cdots,i_q, 1,2,\cdots \hat i_1,\cdots,\hat i_2,\cdots \hat i_2,\cdots,\hat i_q,\cdots,n).$$
Here ``hat'' above symbol means the deletion of this symbol.

The operators $R^q$'s satisfy
$$
\leqalignno{ R^q\cdot R^{n-q}= (-1)^{q(n-q)}Id.&&(3.7)\cr}$$ 

With the help of the operators $R^q$ of an oriented Riemannian manifold
of dimension $n,$ one defines the fiberwise scalar product
 $\Omega(M)^q\times \Omega^q(M) \to \Omega^0(M)$ and the formal adjoints 
$$\delta^{q+1}, \delta^{q+1}(t) :\Omega^{q+1}(M) \to \Omega^q(M),$$
$$({\iota_X^{q-1})}^{\sharp}:\Omega^{q-1}(M)\to \Omega^q(M), \text{and}  
\ {(L_X^q)}^{\sharp}:\Omega^{q}(M)\to \Omega^q(M)$$ of $d^q, d^q(t), \iota_X^q, L_X^q$
by: 

$$\leqalignno{ \ll \omega_1, \omega_2 \gg= (R^n)^{-1}(\omega_1\wedge R^q
(\omega_2)),&&(3.8)\cr}$$

$$\leqalignno{ 
\delta^{q+1} & = (-1)^{nq+1} R^{n-q}\cdot d^{n-q-1}\cdot R^{q+1},\  &(3.9)\cr 
\delta^{q+1} &(t) = (-1)^{nq+1} R^{n-q}\cdot d^{n-q-1}(t)\cdot R^{q+1}, \cr
(\iota_X^q)^{\sharp}&= (-1)^{nq-1} R^{n-q}\cdot \iota_X^{n-q-1}\cdot R^{q-1},\cr
(L_X^q)^{\sharp}&= (-1)^{(n+1)q+1} R^{n-q}\cdot L_X^{n-q}\cdot R^{q}  \cr}
$$

These operators  satisfy:

$$\leqalignno{ 
\ll d\omega_1, \omega_2\gg & = \ll \omega_1, \delta \omega_2\gg,\  &(3.10)\cr 
\ll d(t)\omega_1, \omega_2\gg & = \ll \omega_1, \delta (t)\omega_2\gg, \cr
\ll \iota_X \omega_1, \omega_2\gg & = \ll \omega_1, (\iota_X)^{\sharp} 
\omega_2 \gg,\cr
\ll L_X \omega_1, \omega_2\gg & = \ll \omega_1, (L_X)^{\sharp} \omega_2\gg
,}$$
and
$$\leqalignno{ (L_X)^{\sharp} = (\iota_X )^{\sharp}\cdot \delta + \delta \cdot 
(\iota_X)^{\sharp}.
&&(3.11)\cr}$$

Note that
$ L_X^q + (L_X^q)^{\sharp}$ is a zeroth order differential operator. 
Let $X^{\sharp}$ denote the element in $\Omega^1(M)$ defined by 
$X^{\sharp}(Y):= \ll X, Y \gg\ $ 
and for  $\omega \in \Omega^1(M)$ 
let $E_{\omega}^q: \Omega^q(M) \to \Omega^{q+1}(M),$ denote
the exterior product by $\omega.$ Then we have 
$$\leqalignno{(\iota_X^q)^{\sharp}= \ E_{X^{\sharp}}^{q-1}. 
&&(3.12)\cr}$$

It is easy to see that the scalar products $\ll.,.\gg$ and the operators 
$\delta^q ,\delta^q(t), \iota_X^{\sharp}$ and $ L_X^{\sharp}$ are independent 
of 
the orientation of $M.$  Therefore they are defined (first locally and then being differential operators globally) for an arbitrary 
Riemannian 
manifold, not necessary 
orientable, and satisfy  
 (3.8), (3.10)-(3.12) above.
 
For a Riemannian manifold $(M,g)$ one introduces the scalar product 
$\Omega^q(M)\times \Omega^q(M)\to \Bbb C$
by 

$$\leqalignno{< \omega, \omega'> := \int_M \omega \wedge \omega' = 
\int_M \ll \omega, \omega'\gg dvol(g). &&(3.13)\cr}$$

In view of (3.10), $\delta^{q+1}(t), (\iota_X^q)^\sharp$ and $(L_X^q)^\sharp$ 
are 
formal adjoints of  $d^q(t), \ \iota_X^q(t)$ and $L_X^q$ with respect to the 
scalar product $<.,.>.$ 

For a Riemannian manifold $(M,g),$ one introduces the second order differential operators 
$\Delta_q: \Omega^q(M) \to \Omega^{q}(M),$ the Laplace Beltrami operator,
and $\Delta_q(t): \Omega^q(M) \to \Omega^{q}(M),$ the Witten Laplacian for the function
$h,$ by 
$$\Delta_q:= \delta^{q-1} \cdot d^q + d^{q-1}\cdot \delta^q,$$ and 
 $$\Delta_q(t):= \delta^{q-1}(t)\cdot d^q(t) + d^{q-1}(t)\cdot \delta^q(t).$$ 
Note that $\Delta_q(0)= \Delta_q.$ In view of 
 (3.1) -(3.8) and (3.10) one verifies

$$\leqalignno{\Delta_q(t)= \Delta_q + t(L_{-grad_gh} +L_ {-grad_gh}^{\sharp}) + t^2 ||grad_gh|| Id
&&(3.14)\cr}$$
and that $L_{-grad_gh} +L_ {-grad_gh}^{\sharp}$ is a zeroth  order differential operator.

The operators $\Delta_q(t)$ are  elliptic selfadjoint
and positive, hence their spectra $\text{spect} \Delta_q(t),$ lie on 
$[0,\infty).$ Further, as 
$$ \ker \Delta_q(t)= \{ \omega \in \Omega^q(M) | d^q(t)=0, \delta^q(t)=0 \}$$
one can see that for all $t\geq 0$ $\ker \Delta_q(t)$ is isomorphic 
to $\ker \Delta_q(0).$ 
Hence if $0$ is an eigenvalue of $\Delta_q(0),$ then it is an eigenvalue of 
$\Delta_q(t)$ for all $t$ and with the same 
multiplicity.

A very important fact in the proof of Theorems 3.1 and 3.2 below is 
is that 
$\Delta_q(t) -\Delta_q$ is a zeroth order operator for any $t$.

The following result 
is essentially due to E.Witten, and provides the first main result of the
WHS-theory. 

\proclaim {Theorem 3.1} Suppose that $\tau=(g,h)$ is a 
generalized triangulation of the closed Riemannian manifold M. There exist the 
constants $C_1, C_2, C_3$ and $T_0$ depending on $\tau,$  so that for any 
$t> T_0,$ 
$spect \Delta_q(t) \subset [0,C_1 e^{-C_2t}] \cup [C_3 t,\infty)$ and
the number of the eigenvalues of $\Delta_q(t)$
in the interval $[0,C_1 e^{-C_2t}]$ counted with their multiplicity is 
equal to the number of critical points of index $q.$
\endproclaim

The above theorem states the existence of a gap in the spectrum of 
$\Delta_q(t),$ namely the open interval $(C_1 e^{-C_2t}, C_3t) ,$
which widens to $(0,\infty)$ when $t\to \infty.$

Clearly $C_1, C_2, C_3$ and  $T_0$ determine a constant $T$, so that 
$1\in (C_1e^{-C_2t}, C_3t)$
and for $t\geq T,$  
$$spect \Delta_q(t) \cap [0, C_1 e^{-C_2 t}] = spect \Delta_q(t) 
\cap [0, 1] $$ and 
$$spect \Delta_q(t) \cap [C_3t, \infty) = spect \Delta_q(t) \cap [1,\infty).$$

For $t>T$ we denote by $\Omega^q(M)(t)_{sm}$ the finite dimensional subspace
of dimension $m_q,$ the number of critical points of index $q,$
generated by the $q-$eigenforms 
of $\Delta_q(t)$ corresponding to the eigenvalues of $\Delta_q(t)$ 
smaller than $1.$
The elliptic theory implies that these eigenvectors, a priori elements in the 
$L_2-$completion of $\Omega^q(M),$ are actually in $\Omega^q(M).$ Note that 
$d(t)(\Omega^q(M)(t)_{sm})\subset \Omega^{q+1}(M)(t)_{sm},$
so that $(\Omega^*(M)(t)_{sm}, d^*(t))$ is a finite dimensional cochain subcomplex
of $(\Omega^*, d^*(t))$ and $e^{th} ( \Omega^*(M)(t)_{sm}, d^*(t))$ is a 
finite dimensional subcomplex of $(\Omega^*(M), d^*).$

For $t>T,$ consider the composition of morphisms of cochain complexes  denoted by 
$l^*(t),$
$$ (\Omega^*(M)_{sm}, \overline{d}^*(t))@>S^*(t)>> (\Omega^*(M)_{sm}, d^*(t))
@>e^{th}>> (\Omega^*(M), d^*) @> Int^*>>(C^*(M,\tau),\partial^*),$$
with $S^q(t)= (\frac{\pi}{t})^{\frac{n-2q}{4}} e^{tq}Id,$ and 
$\overline{d}^q(t): =  (\frac{\pi}{t})^{1/2} e^{-t} d^q(t).$
$S^*(t)$ is an isomorphism of cochain complexes referred to as the
"rescaling isomorphism". The following theorem due to Helffer-Sj\"ostrand, 
cf [HS2], provides the second main result of the WHS-theory.

\proclaim {Theorem 3.2} (Helffer-Sj\"ostrand) Given $M$ a 
closed manifold and $\tau=(g,h)$ a generalized triangulation, there 
exists $T_1 >0,$ depending
on $\tau,$ so that for $t> T_1$
$1\notin \text{spect} \Delta_q(t)$ and $l^*(t)$ is an isomorphism of cochain 
complexes. 

Moreover, for $t>T_1$ there exists a family of isometries 
$\Cal J^q(t): C^q(M,\tau) \to \Omega^q(M)(t)_{sm}$ of finite dimensional vector spaces so that 
$l^q(t) \Cal J^q(t) = Id + O(1/t).$ It is understood that 
$C^q(M,\tau)$ is equipped with 
the canonical scalar product defined in section 1, before Theorem 1.2,
and $\Omega^q(M)(t)_{sm}$ with the scalar product $< .,.>$ defined 
by (3.13).
\endproclaim

Theorem 3.2 provides inside $ (\Omega^*(M), d^*(t)),$ (a reparametrization of 
$ (\Omega^*(M), d^*)$
induced by the multiplication operator $e^{th} S^q(t): 
\Omega^q(M) \to \Omega^q(M)$) the 
finite dimensional subcomplex  $(\Omega^*(M)_{sm}, \overline{d}^*(t))$ which, 
after rescaling, is asymptotically 
isometric to $(C^*(M,\tau),\partial^*).$

Recall that De Rham Hodge theory provides a canonical and unique 
representation of each cohomology class of $ (C^*(M,\tau), \partial^*)$ by 
harmonic $q-$forms with respect to $g$. Theorem 3.2 provides, asymptotically, 
a  canonical and unique
representation  of the full complex $(C^*(M,\tau),\partial^*)$ 
and its base $E_x$ inside  $(\Omega^*(M), d^*).$

\vfill \eject 

\proclaim{4. Ideas of the proof of Theorems 3.1 and 3.2}  \endproclaim 

The proof of Theorems 3.1 and 3.2 is based on a mini-max criterion for detecting a gap
in the spectrum of a positive selfadjoint operator in a Hilbert space $H,$
 Lemma 4.1 below, and on the  explicit formula for $\Delta_q(t)$ in admissible 
 coordinates in a neighborhood of the critical points. 

 \proclaim {Lemma 4.1} Let $A: H \to H$ be a densely defined (not necessary bounded ) 
 self adjoint positive operator in a Hilbert space $(H,<,>)$ 
and  $a,b$ 
two real numbers so that $0<a <b < \infty.$  Suppose that there exists 
two closed subspaces $H_1$ and $H_2$ of $H$ with $H_1\cap H_2 =0$ 
and $H_1 + H_2 =H$ such that:
 \newline (1) $<Ax_1, x_2> \ \leq \ a||x_1||^2$ for $x_1\in H_1,$
 \newline (2) $<Ax_1, x_2> \ \geq \ b||x_2||^2$ for $x_1\in H_2.$

 Then $spect A\bigcap (a,b)= \emptyset.$

 \endproclaim

 The proof of this Lemma is elementary and is left as an exercise for the reader.

Consider $x\in Cr(h)$ and choose admissible coordinates $ (x_1,x_2, ...,x_n) $
in the neighborhood of $x$. Since with respect to these coordinates 
$$h(x_1,x_2, ...,x_n) = k  - 1/2 (x_1^2+\cdots+ x_k^2) +1/2 (x_{k+1}^2+\cdots + x_n^2)$$
 and  $g_{ij}(x_1,x_2, ...,x_n) = \delta_{ij},$ by (3.14) 
the operator $\Delta_q(t)$
 has the form: 
$$\leqalignno{\Delta_{q,k}(t)= \Delta_q + t M_{q,k} +t^2 (x_1^2 +\cdots +x_n^2)Id
&&(4.1)\cr}$$
with $$\Delta_q (\sum _ I a_I(x_1,x_2, ...,x_n) dx_I)
= - (\sum _{i=1}^n \frac{\partial^2}{\partial x_i^2} a_I (x_1,x_2, ...,x_n)) dx_I,$$
and $M_{q,k}$ the linear operator determined by 
$$
\leqalignno{ M_{q,k} (\sum _ I a_I(x_1,x_2, ...,x_n) dx_I)= 
\sum _I \epsilon_I^{q,k} a_I(x_1,x_2, ...,x_n) dx_I .&&(4.2)\cr}$$
Here  $I =(i_1,i_2\cdots i_q),$  $1\leq i_1 < i_2 \cdots< i_q \leq n,$
 $ d_I= dx_{i_1}\wedge \cdots \wedge d_{i_q} $ and 
$$\epsilon _I^{q,k} = -n+2k-2q +4 \sharp\{j|k+1 \leq i_j \leq n\},$$ 
where $\sharp A$ denotes the  cardinality of the set $A.$
Note that $\epsilon^{q,k}_I \geq -n$ and is $=-n$ iff $q=k.$

Let $\Cal S^q(\Bbb R^n)$ denote the space 
of smooth $q-$forms  
$\omega = \sum_I a_I(x_1,x_2, ...,x_n) dx_I$ 
with $a_I(x_1,x_2, ...,x_n)$ rapidly
decaying functions. The operator $\Delta_{q,k}(t)$ acting on 
$\Cal S^q(\Bbb R^n)$ is globally elliptic (in the sense of [Sh1] 
or [H\"o]), selfadjoint and positive. This operator is the harmonic 
oscillator in $n$ variables acting on $q-$forms and its properties can 
be derived 
from the harmonic oscillator in one variable $-\frac{d^2}{dx^2}+a +bx^2$
acting on functions. In particular the following result holds.

\proclaim{ Proposition 4.2} (1) $\Delta_{q,k} (t),$ regarded as an unbounded 
densely defined operator on the $L_2-$completion of  $\Cal S^q(\Bbb R^n),$ is 
selfadjoint, positive and its spectrum is contained in $2t\Bbb Z_{\geq 0}$
(i.e positive integer multiple of $2t$).

(2) $\ker \Delta_{q,k}(t)=0\  \text{if} \ 
k\ne q$ and $dim \ker\Delta_{q,q}(t) = 1.$

(3) $\omega_{q,t}= (t/{\pi})^{n/2} e^{-t\sum_i x_i^2/2} dx_1\wedge\cdots \wedge dx_q$ 
is the generator of $\ker\Delta_{q,q}(t)$ with the $L_2-$norm $1.$
\endproclaim

For details consult [BFKM] page 805.

Choose a smooth  function $\gamma_\eta(u),$ 
$\eta\in (0,\infty), \ u\in \Bbb R,$
which satisfies :
$$ \leqalignno{\gamma_\eta(u)=\quad  
\left\{\aligned 1\ \text {if} \ & u\leq \eta/2 \\ 0 \ \text{if}\  &u 
>\eta\endaligned \right\}.&&(4.3)\cr}$$

Introduce 
$\tilde \omega_{q,t}^\eta \in \Omega^q_c(\Bbb R^n) $defined by 

$\tilde\omega_{q,t}^\eta (x)= \beta_q^{-1}(t) \ \gamma_\eta (|x|) 
\omega_{q,t}(x)$ with $|x|= \sqrt{\sum_i x_i^2}$ and

$$\leqalignno{ \beta_q(t)= (t/{\pi})^{n/4}( \int_{\Bbb R^n} 
\gamma_\eta^2(|x|) e^{-t \sum_ix_i^2} dx_1\cdots dx_n)^{1/2}.
&&(4.4)\cr}$$

The smooth form $\tilde\omega_{q,t}^\eta$ has the support in the disc 
of radius $\eta,$ agrees with $\omega_{q,t}$ on the disc of 
radius $\eta/2$ and satisfies 

$$\leqalignno{< \tilde\omega_q^\eta(t), \tilde\omega_q^\eta(t)> =1
&&(4.5)\cr}$$ with respect to the scalar product $<.,.>$ on $\Cal S^q(\Bbb R^n)$
induced by the Euclidean metric. The following proposition can be obtained 
by elementary calculations in coordinates in view of the explicit formula of
$\Delta_{q,k}(t)$ cf [BFKM], Appendix 2.

\proclaim{Proposition 4.3} For a fixed $r\in \Bbb N_{\geq 0}$ there
exists $C,C',C'', T_0, \epsilon_0$ so that $t >T_0$ and 
$\epsilon <\epsilon_0$ imply 

(1) $ |\frac{\partial^{|\alpha|}}{\partial x_1^{\alpha_1}\cdots 
\partial x_n^{\alpha_n}}\Delta_{q,q}(t) 
\tilde\omega_{q,t}^\epsilon(x)|\leq Ce^{-C't}$ for any 
$x\in \Bbb R^n$ and multiindex $\alpha= (\alpha_1,\cdots,\alpha_n)$, with 
$|\alpha|=\alpha_1+\cdots+ \alpha_n \leq  r.$

(2) $<\Delta_{q,k}(t)\tilde\omega_{q,t}^{\epsilon}, 
\tilde \omega_{q,t}^{\epsilon}>\  \geq\  2t|q-k|$

(3) If $\omega \perp \tilde\omega^\epsilon _{q,t}$  with respect to the scalar product $<.,.>$
then $$<\Delta_{q,q}\omega,\omega>\  \geq \ C'' t||\omega||^2.$$

\endproclaim

For the proof of Theorems 3.1 and 3.2 we choose 
$\epsilon_0 >0$ so that for each $y\in Cr(h)$ there exists an 
admissible coordinate chart
$\varphi_y: (U_y, y) \to (D_{2\epsilon},0) $ so that
$U_y\cap U_z=\emptyset$ for $y\ne z,$ $y,z\in Cr(h).$

Choose once for all such an admissible coordinate chart for any $y\in Cr_q(h).$
Introduce the smooth forms $\overline{\omega}_{y,t}\in \Omega^q(M)$ 
defined by 
$$\leqalignno{ \overline{\omega}_{y,t} |_{M\setminus \varphi_y^{-1}
(D_{2\epsilon})}=0, \ 
\overline{\omega}_{y,t}|_{\varphi_y^{-1}(D_{2\epsilon})}= 
\varphi_y^*(\tilde\omega^\epsilon_{q,t}).
&&(4.6)\cr}$$

The forms $\overline{\omega}_{y,t}\in\Omega^q(M), \ y\in Cr_q(h)$ are 
orthonormal. Indeed 
if $y\ne z, \  y,z\in Cr_q(h),$ $\overline{\omega}_{y,t}$ and 
$\overline{\omega}_{z,t}$ 
have disjoint support, hence are orthogonal, and because the support 
of $\overline{\omega}_{y,t}$
is contained in an admissible chart, 
$< \overline{\omega}_{y,t}, \overline{\omega}_{y,t}>=1$ by (4.5).

For $t>T_0,$ with $T_0,$ given by Proposition 4.3, we define 
$J^q(t): C^q(X,\tau) \to \Omega^q(M)$
to be the linear map determined by 
$$J^q(t)(E_y) = \overline{\omega}_y(t),$$ where 
$E_y\in C^q(X,\tau)$ is given by $E_y(z)= \delta_{yz}$ for $y,z\in Cr(h)_q.$ 
$J_q(t)$ is an isometry, thus in particular injective.

{\bf Proof of Theorems 3.1 and 3.2:} (sketch). Take $H$ to be the 
$L_2-$completion of 
$\Omega^q(M)$
with respect to the scalar product $<., .>,$ $H_1:= J^q(t)(C^q(M,\tau))$ and $H_2= H_1^{\perp}.$
Let $T_0, C, C', C''$ be given by Proposition 4.3 and define 
$$C_1:= \inf _{z\in M'} || grad_g h(z)||,$$ with 
$\ M'= M\setminus \bigcup_{y\in Cr_q(h)} \varphi_y^{-1}(D_\epsilon),$
$$C_2= \sup_{x\in M} ||(L_{ -grad_g h} +L_ {-grad_g h}^{\sharp})(z)||;$$  
here $||grad_g h(z) ||$ resp.
$||(L_{-grad_g h} +L_ {-grad_g h}^{\sharp})(z)||$ denotes the norm of
the vector $ grad_g h(z)\in T_z(M)$ resp. of the linear map 
$(L_{-grad_g h} +L_ {-grad_g h}^{\sharp})(z): \Lambda^q(T_z(M)) \to \Lambda^q(T_z(M))$ 
with respect to the scalar product  induced in $T_z(M)$ and 
$\Lambda^q(T_z(M))$ by $g(z).$
Recall that if $X$ is a vector field then $L_X +L_X^{\sharp}$ is a zeroth 
order differential operator,
hence an endomorphism of the bundle $\Lambda^q(T^*M)\to M.$

We can use the constants 
$T_0, C, C', C'', C_1, C_2$ to construct $C'''$ and $\epsilon_1$ so that for $t> T_0$ and $\epsilon <\epsilon_1,$  
we have 
$< \Delta_q(t) \omega, \omega > \geq C_3t <\omega, \omega>$ for any $\omega \in H_2$
 (cf. [BFKM], page 808-810).
 
Now one can apply Lemma 4.1 whose hypotheses are satisfied for 
$a= Ce^{-C't}, b =C''' t$ and $t>T_0.$
This concludes the first  part of Theorem 3.1. 

 Let $Q_q(t),$  $t>T_0$ denote the orthogonal projection in $H$ on the span of  the eigenvectors 
 corresponding the eigenvalues smaller than $1.$ In view of the ellipticity of $\Delta_q(t)$ all 
 these eigenvectors are smooth $q-$forms. An additional important estimate is 
given by the following Proposition:

\proclaim{Proposition 4.4} For $r\in \Bbb N_{\geq 0}$ one can find 
$\epsilon_0 >0 $ and $C_3, C_4$ so that 
for $t> T_0$ as constructed above, and any $\epsilon <\epsilon_0$
one has, for any $v\in C^q(M,\tau)$
$$ \sup_{x\in M}||(Q_q(t) J^q(t) - J^q(t)) (v)|| \leq C_3 e^{- C_4t} 
||v||,$$
$(Q_q(t) J^q(t)- J^q(t))(v) \in \Omega^q(M),$
with similar estimates for the 
$C^p-$norm of $(Q_q(t) J^q(t) - J^q(t),$ with $p\leq r$.
\endproclaim

The proof of this Proposition is contained in [BZ1], page 128 and [BFKM] page 811. Its proof requires (3.14), Proposition 4.3 and general estimates coming from the ellipticity of $\Delta_q(t).$  

Proposition 4.4 implies that for $t$ large enough, say $t> t_0,$
$\Cal I^q(t) := Q_q(t) J^q(t)$ is bijective, which  finishes the proof of 
Theorem 3.1.

For $t\geq T_0,$ as constructed in Proposition 4.4, let $\Cal J^q(t)$ be 
the isometry defined by 
$$\leqalignno{ \Cal J^q(t):= \Cal I^q(t) (\Cal I^q(t)^{\sharp} \Cal I^q(t))^{-1/2}
&&(4.7)\cr}$$
and denote by $U_{y,t}:= \Cal J^q(t) (E_y) \in \Omega^q(M), y\in Cr(h)_q.$  
Proposition 4.4 implies that there 
exists $t_0$ and 
$C$ so that for $t > t_0$  and $y\in Cr(h)_q$ one has:

$$\leqalignno{ \sup_{z\in M\setminus \varphi_y^{-1}(D_\epsilon)} 
||U_{y,t}(z)|| \leq Ce^{-\epsilon t}, &&(4.8)\cr}$$

$$\leqalignno{ ||U_{y,t}(z) - \overline\omega_{y,t}(z)|| 
\leq C \frac{1}{t},\  \text{for} \  z\in W_y^-\cup \varphi_y^{-1}(D_\epsilon). &&(4.9)\cr}$$

To check Theorem 3.2 it suffices to show that 

$$| \int_{W_{x'}^-} U_{x,t} e^{th} -(\frac{t}{\pi})^{\frac{n-2q}{4}} e^{tq}
\delta_{xx'} | 
\leq C''\ \frac{1}{t} $$
for some $C''> 0$ and any $x,x'\in  Cr(h)_q.$

If $x\ne x'$ this follows from (1). If $x=x'$ from (4.8) and (4.9).
\hfill $\square$

\vfill \eject 

\proclaim{5. Extensions and a survey of applications  }  \endproclaim 

1. One can relax the definition of the generalized triangulation by dropping 
from 
$C1$ the constraint on $g$ to be Euclidean in the neighborhood of the critical 
points. 
This will keep Theorems 3.1 and 3.2 valid as stated; however almost all calculations 
will be longer since the explicit formulae for $\Delta_q(t)$ and its spectrum 
when regarded on $\Cal S^*(\Bbb R^n)$ will be more complicated.

2.One can drop the hypothesis that the Morse function $h$ is self indexing.
In this case  Theorem
3.1 remains true as stated but 
in Theorem 3.2, 
$l^q(t)$ should be replaced by 

$$  (\Omega^*(M)_{sm}, d^*(t))
@>e^{th}>> (\Omega^*(M), d^*) @> Int^*>>(C^*(M,\tau),\partial^*),$$

and $\Cal J^q(t)$ by $\Cal J^q(t)\cdot \Sigma^q(t)$ with 
$\Sigma^q(t): (C^q(M,\tau),\partial^*)\to (C^q(M,\tau), \partial^*(t)),$ and 
$\partial^*(t)= \Sigma^{q+1}(t)\cdot \partial^q \cdot (\Sigma^q(t))^{-1}$ 
the morphism of cochain complexes defined by          
$$ \Sigma^q(t) (E_x) = (\frac{\pi}{t})^{\frac{n-2q}{4}} e^{ h(x)}E_x,$$
$x\in Cr_q(h)$ 
cf. [BFK3].

3. One can twist both complexes $(C^*(M,\tau), \partial^*)$ and 
$(\Omega^*, d^*)$
by a finite dimensional representation of the fundamental group, 
$\rho: \pi_1(M) \to GL(V).$
In this case an additional data is necessary: a Hermitian structure $\mu$ 
on the flat bundle $\xi_\rho$ induced by 
$\rho.$
The "canonical" scalar product on $(C^*(M,\tau,\rho), \partial^*)$
will be obtained by using the critical points (the cells of the 
generalized triangulation 
and the hermitian scalar product provided by $\mu$ in the fibers
of $\xi_\rho$ above the critical points. The  De-Rham complex in this case is replaced by
$(\Omega^*(M,\rho), d_\rho^*)$ of differential forms with coefficients in 
$\xi_\rho$ and the differential is provided by the flat connection in 
$\xi_\rho.$ The scalar product will require in addition of the Riemannian metric 
$g$ the Hermitian structure $\mu$ (cf. [BFK] or [BFK2]).
Under the hypotheses that the Hermitian structure is parallel in small 
neighborhoods of the critical points, the  proofs of Theorems 3.1 and 3.2 
remain the same.
An easy continuity argument permits to reduce the case of an arbitrary 
Hermitian structure
to the previous one by taking a $C^0$ approximation of a given Hermitian 
structure by Hermitian structures which are parallel near the critical points.
Since
the Witten Laplacians do not involve derivatives of the Hermitian structure
such a reduction is possible.
If the representation is a unitary representation in a finite dimensional Euclidean space 
one has a canonical Hermitian structure in $\xi_\rho$ (parallel with respect to the flat 
canonical connection in $\xi_\rho.$)
This    
   extension was used in the new proofs of the Cheeger- Muller theorem 
and its extension
concerning the  comparison of the analytic and the Reidemeister torsion. cf. [BZ], [BFK].

4.One can further extend the WHS-theory to the case where  $\rho$ is a special 
type of an infinite 
dimensional representation, a representation of the fundamental group in an $\Cal A-$
Hilbert module of finite type. This extension was done in [BFKM] for $\rho$
unitary and in [BFK4] for
$\rho$ arbitrary. In this case the Laplacian $\Delta_q(t)$ do 
not have discrete spectrum 
and it seems quite remarkable that Theorems 3.1 and 3.2 remain true. It is even more 
surprising that 
exactly the same arguments as presented above can be adapted to prove them.
A particularly interesting situation is the case of the left regular representation of a countable group
$\Gamma$ on the Hilbert space $L_2(\Gamma)$ when regarded as an 
$\Cal N(\Gamma)$ right 
Hilbert module of the von Neumann algebra $\Cal N(\Gamma),$  
cf.[BFKM] for definitions.
One can prove that Farber extended $L_2-$cohomology of $M,$ a compact smooth
manifold with infinite fundamental group defined analytically (i.e. using differential forms and a Riemannian metric) and combinatorially (i.e using a triangulation) are isomorphic and therefore 
the classical
$L_2-$Betti numbers and Novikov-Shoubin invariants defined analytically and combinatorially are the same. For the last fact see
see [BFKM].

This WHS-theory was a fundamental tool in the proof of the equality of the
\newline $L_2-$analytic and the $L_2-$Reidemeister 
torsion presented [BFKM].

5. One can further extend Theorems 3.1 and 3.2 to bordisms $(M,\partial_-M, 
\partial_+M),$
and $\rho$ a representation of $\Gamma=\pi_1(M)$ on an $\Cal A-$Hilbert module of finite type.
In this case one has first to extend the concept of generalized triangulation  to such bordisms. This will 
involve a pair $(h, g)$ which in addition to the requirements C1-C3 is 
supposed to satisfy the following assumptions:
$g$ is product like near $\partial M=\partial_-M \cup \partial_+M,$

$h:M \to [a, b]$ with $h^{-1}(a)= \partial_-M,$  $h^{-1}(b)= \partial_+M,
\  a,b$ regular values, and $h$ linear on the geodesics normal to 
$\partial M$ near $\partial W.$
This extension was done in [BFK2] and was used to prove gluing 
formulae for 
analytic torsion and to extend the results of [BFKM] to manifolds 
with boundary.

5. One can actually extend WHS-theory to the case where $h$ is 
generalized Morse function, i.e. the critical points are either 
nondegenerated or birth-death. This extension is much more subtle and 
very important. Beginning work in this direction was done by Hon Kit Wai 
in his OSU dissertation.
\vfill \eject 

\proclaim{ 6. References} \endproclaim

\end

\enddocument